\documentclass{amsart}[11pt]

\newtheorem{theorem}{Theorem}[section]
\newtheorem{lem}[theorem]{Lemma}
\newtheorem{prop}[theorem]{Proposition}
\newtheorem{cor}[theorem]{Corollary}
\newtheorem*{thA}{Theorem A}
\newtheorem*{thB}{Theorem B}
\newtheorem*{thC}{Theorem C}

\theoremstyle{definition}
\newtheorem{definition}[theorem]{Definition}

\newtheorem{claim}{Claim}

\theoremstyle{remark}

\numberwithin{equation}{section}

\newcommand{\half}{\frac{1}{2}}

\newcommand{\Hor}{{\mathcal{H}}}
\newcommand{\V}{{\mathcal{V}}}

\newcommand{\nb}{\nu(\Sigma)}
\newcommand{\nbp}{\nu_p(\Sigma)}

\newcommand{\ra}{\rightarrow}
\newcommand{\R}{{\mathbb{R}}}

\newcommand{\Soul}{\Sigma}

\newcommand{\Ddt}{\frac{\text{D}}{\text{dt}}}

\newcommand{\dz}{\frac{\text{d}}{\text{dt}}\bigg|_{t=0}}

\newcommand{\grad}{\text{grad}}

\newcommand{\hess}{\text{hess}}

\newcommand{\lb}{\langle}
\newcommand{\rb}{\rangle}
\newcommand{\expp}{\exp^{\perp}}

\newcommand{\RN}{R^{\nabla}}

\newcommand{\bX}{\bar{X}}

\newcommand{\bU}{\bar{U}}
\newcommand{\bV}{\bar{V}}
\newcommand{\bW}{\bar{W}}
\newcommand{\RC}{\text{\r{R}}}

\begin{document}

\title[Conditions for Nonnegative Curvature]{Conditions for Nonnegative
Curvature on Vector Bundles and Sphere Bundles}
\author{Kristopher Tapp}
\address{Department of Mathematics\\ SUNY at Stony Brook\\
         Stony Brook, NY 11794-3651}
\email{ktapp@math.sunysb.edu}

\subjclass{Primary 53C20}
\date{\today}
\keywords{nonnegative curvature, soul, connection, vector bundle}

\begin{abstract}
This paper addresses Cheeger and Gromoll's question of which vector bundles
admit a complete metric of nonnegative curvature, and relates their question to the issue
of which sphere bundles admit a metric of positive curvature.  We show that
any vector bundle which admits a metric of nonnegative curvature must admit a connection,
a tensor, and a metric on the base space which together satisfy a certain differential
inequality.  On the other hand, a slight sharpening of this condition is sufficient
for the associated sphere bundle to admit a metric of positive curvature.
Our results sharpen and generalize Walschap and Strake's conditions under which a
vector bundle admits a connection metric of nonnegative curvature.
\end{abstract}

\maketitle

\section{Introduction}
A well-known question in Riemannian geometry is to what extent the converse of Cheeger
and Gromoll's soul theorem holds.
Their theorem states that any complete noncompact manifold, $M$,
with nonnegative sectional curvature is diffeomorphic to the normal bundle of a compact
totally geodesic submanifold, $\Soul\subset M$, called the ``soul of $M$''~\cite{CG}.
The converse
question is the classification problem:  which vector bundles over compact nonnegatively
curved base spaces can admit complete metrics of nonnegative curvature?
There are vector bundles which are know not to admit nonnegative curvature, but
in all such examples the base space 
has an infinite fundamental group~\cite{WO},\cite{Tapp3},\cite{BK}.
Trivial positive results include all vector bundles over $S^1$, $S^2$ and $S^3$,
$TS^n$ for any $n$, and more generally all homogeneous vector bundles over homogeneous spaces.
As for nontrivial positive results,
Yang obtained nonnegatively curved
metrics on rank 2 vector bundles over $CP^n\#\overline{CP^n}$~\cite{Yang}. 
More recently, Ziller and Grove constructed nonnegatively curved metrics on all vector bundles
over $S^4$ and $S^5$~\cite{GZ}.

Our first result is a necessary condition for a vector bundle to admit a metric of
nonnegative curvature.
Suppose that $M$ is an open (ie., complete and noncompact)
manifold with nonnegative curvature, and $\Soul$ is a soul of $M$.
Let $\nb$ denote the normal bundle of $\Soul$, let $p\in\Soul$, $X,Y\in T_p\Soul$,
and $W,V\in \nbp$.  Let $\nabla$ denote the connection in $\nb$, and let $\RN$ denote its
curvature tensor, so that $\RN(X,Y)W\in\nbp$.  We can enlarge
the domain of the tensor $\RN$ by defining $\RN(W,V)X$ to be the vector
in $T_p\Soul$ for which $\lb\RN(W,V)X,Y\rb=\lb\RN(X,Y)W,V\rb$ for every $Y\in T_p\Soul$.
Let $k_{\Soul}$ describe the unnormalized sectional
curvatures of $\Soul$; that is, $k_{\Soul}(X,Y)=\lb R(X,Y)Y,X\rb$.
Similarly, let $k_F(W,V)=\lb R(W,V)V,W\rb$ describe unnormalized sectional curvatures of
2-planes perpendicular to $\Soul$ (``F'' stands for ``fiber'', since $k_F$
really describes the curvatures of the fibers at points of $\Soul$).
By parallel transporting $W$ and $V$ along geodesics
from $p$ in $\Soul$, we can regard $k_F(W,V)$ as a real valued function on $\Soul$ near $p$;
by $\hess_{k_F(W,V)}(X)$ we denote the hessian of this function
in the direction $X$.  We think of \{$\RN$, $k_{\Soul}$, $k_F$\} as the structure of $M$
which is visible at points of the soul.  Our necessary condition for nonnegative curvature
is the following relationship between these visible structures:

\begin{thA}
If $M$ is an open manifold of nonnegative curvature with soul $\Soul$, then for
any $p\in\Soul$, $X,Y\in T_p\Soul$, and $W,V\in\nbp$,
$$
\lb (D_X\RN)(X,Y)W,V\rb^2\leq(|\RN(W,V)X|^2+(2/3)\hess_{k_F(W,V)}(X))\cdot k_{\Soul}(X,Y).
$$
\end{thA}

An obvious question is whether the condition in Theorem A is sufficient; that is,
if a vector bundles admits structures statisfying the inequality of the theorem,
then must it admit a metric of nonnegative curvature?  We discuss first the case
of connection metrics, about which much is allready known.

A ``connection metric'' $g_E$ on the total space
$E$ of a vector bundle $\R^k\ra E\stackrel{\pi}{\ra}\Soul$ is a metric arising
from the following construction.  Choose a Euclidean structure $\lb\cdot,\cdot\rb$ on
the bundle (which means a smoothly varying choice of inner products on the fibers),
a connection $\nabla$ that is compatible with the Euclidean structure, a metric
$g_{\Soul}$ on $\Soul$, and a rotationally-symmetric metric $g_0$ on $\R^k$.  Then
there is a unique metric $g_E$ on $E$ for which $\pi:(E,g_E)\ra(\Soul,g_{\Soul})$
is a Riemannian submersion with horizontal distribution, $\Hor$,
determined by $\nabla$, and with totally geodesic fibers isometric to $(\R^k,g_0)$.
By a connection metric $g_{E^1}$ on the total space $E^1$ of the associated sphere budle
$S^{k-1}\ra E^1\stackrel{\pi}{\ra}\Soul$, we mean the intrinsic metric induced on the sphere
of radius 1 about $\Soul$ in $(E,g_E)$.

In Theorem A, if the metric on $M$ is a connection metric, then $k_F$ is parallel, so
$\hess_{k_F(W,V)}(X)=0$.  The inequality therefore becomes:
\begin{equation}\label{E:ConnectInequality}
\lb (D_X\RN)(X,Y)W,V\rb^2\leq|\RN(W,V)X|^2\cdot k_{\Soul}(X,Y).
\end{equation}
This special case of Theorem A is not new.
In~\cite{Wal}, Strake and Walschap studied conditions under which a vector bundle admits a
connection metric of nonnegative curvature.  Their necessary condition is
stronger than inequality~\ref{E:ConnectInequality}:
\begin{equation}\label{E:SWInequality}
\lb (D_X\RN)(X,Y)W,V\rb^2\leq|\RN(W,V)X|^2\cdot (k_{\Soul}(X,Y)
		-\frac{3}{4}\epsilon^2|\RN(X,Y)W|^2),
\end{equation}
where $\pi\epsilon>0$ is a bound on the diameters of spheres about the origin in $(\R^k,g_0)$.

We prove the following weak converse to Theorem A:

\begin{thB}
Let $\Soul$ be a compact manifold, let $\R^k\ra E\stackrel{\pi}{\ra}\Soul$
be a vector bundle over $\Soul$, and let $S^{k-1}\ra E^1\stackrel{\pi}{\ra}\Soul$
be the associated sphere bundle.

\begin{enumerate}
\item
$E^1$ admits a connection metric $g_{E^1}$ of positive curvature if and only if
the vector bundle admits a Euclidean structure
$\lb\cdot,\cdot\rb$, a compatible connection $\nabla$, and a metric $g_{\Soul}$ on $\Soul$
such that the following inequality holds for all
$p\in\Soul$, $X,Y\in T_p\Soul$ with $X\wedge Y\neq 0$ , and $W,V\in E_p$ with
$W\wedge V\neq 0$:
$$
\lb (D_X\RN)(X,Y)W,V\rb^2 < |\RN(W,V)X|^2\cdot k_{\Soul}(X,Y)
$$
\item
Further, if the vector bundle admits structures for which this inequality is satisfied,
then $E$ admits a complete connection metric $g_E$ of nonnegative curvature.
\end{enumerate}
\end{thB}

Notice the strict inequality implies that $(\Soul,g_{\Soul})$ has positive curvature.
To prove part 2, we show that $g_0$ can be chosen so that the connection
metric $g_E$ on $E$ determined by the data $\{g_{\Soul},\lb\cdot,\cdot\rb,\nabla,g_0\}$
has nonnegative curvature.  Additionally, the boundary of a small
ball about the soul (zero section) of $(E,g_E)$ has positive intrinsic curvature,
which proves one direction of part 1 of the theorem.

We describe next some ways in which Theorem B overlaps known results related to connection
metrics of nonnegative and positive curvature.

\begin{itemize}

\item
One direction of part 1, namely that positive curvature implies the inequality,
follows from the argument in~\cite{Wal} by which Walschap and Strake established
inequality~\ref{E:SWInequality}.  We will elaborate on this remark in Section 7.

\item
Part 2 is an improvement of Strake and Walschap's sufficient condition for a connection
metric of nonnegative curvature, which
is equivalent to our condition with the right side of the inequality multiplied by
$1/2$.

\item
Part 2 of Theorem B in the case where $k=2$ and the vector bundle is oriented was done by
Strake and Walschap
in~\cite{Wal}.  Part 1 of Theorem B in this case follows from Strake and Walschap's
work, and also appears explicitly in~\cite{CDR}.
In this case, $\RN$ can be identified with the 2-form $\Omega$ on $\Soul$
given by $\Omega(X,Y)=\lb \RN(X,Y)W,JW\rb$, where $|W|=1$ and $J$ denotes the almost
complex structure on $E$.  The inequality of Theorem B becomes:
$(D_X\Omega(X,Y))^2<|i_X\Omega|^2\cdot k_{\Soul}(X,Y)$, where $i_X\Omega = \Omega(X,\cdot)$.
Also, inequality~\ref{E:SWInequality} becomes:
$(D_X\Omega(X,Y))^2\leq|i_X\Omega|^2\cdot(k_{\Soul}(X,Y)-\frac{3}{4}\epsilon^2\Omega(X,Y)^2)$;
Yang proved in~\cite[Lemma 1]{Yang} that this last inequality together with the
inequality $k_{\Soul}(X,Y)\geq\frac{3}{4}\epsilon^2\Omega(X,Y)^2$ provide a
necessary and sufficient condition for the connection metric $g_{\epsilon}$
on $E^1$ with fibers of length $2\pi\epsilon$ to have nonnegative curvature.
Since when $k=2$ the sphere bundle $E^1$ is a principal bundle, nonnegative curvature on the
sphere bundle implies nonnegative curvature on the vector bundle.

\item
The strict inequality of Theorem B implies that $\RN(W,V)X=0$ only when
$X=0$ or $W\wedge V=0$.
This is equivalent to saying that the induced connection in the sphere bundle $E^1$ is
``fat''.  The concept of fatness was introduced by Weinstein in~\cite{Wein}.
Among other restrictions, it implies that $\text{dim}(\Soul)$ is even and is $\geq k$,
with equality possible only if $\text{dim}(\Soul)$ is 2, 4, or 8.  
Derdzinski and Rigas proved in~\cite{DR} that the only $S^3$ bundle
over $S^4$ which admits a fat connection is the Hopf bundle
$S^3\ra S^7\ra S^4$.  This result rules out the possibility of using 
Theorem B to produce metrics of positive curvature on 7 dimensional exotic spheres.
We refer the reader to~\cite{Fatness} for a survey of results related to fatness.
Because of the fatness implication, the strict inequality of Theorem B should probably
be considered much stronger than the non-strict inequality of
equation~\ref{E:ConnectInequality}.  
\end{itemize}

We return to the general problem of finding sufficient conditions for nonnegative curvature
on $E$ and for positive curvature on $E^1$.  The inequality
of Theorem A is a relationship between the different curvatures that are visible at the soul, namely,
the curvatures of 2-planes tangent to $\Soul$ (described by $k_{\Soul}$), the curvature
$\RN$ of the connection in $\nb$, and the curvatures of ``vertical'' 2-planes
(described by $k_F$).  It is useful to write
$k_F(W,V)=\lb R_F(W,V)V,W\rb=R_F(W,V,V,W)$, where $R_F$, which we call the ``vertical
curvature tensor'', is just the restriction of the curvature tensor $R$ of $M$ to vectors in
$\nb$, so that $(R_F)_p:(\nbp)^4\ra\R$.

More generally, a tensor $R_F$ on a vector bundle $\R^k\ra E\stackrel{\pi}{\ra}\Soul$
such that for each $p\in\Soul$, the map $(R_F)_p:(E_p)^4\ra\R$ has the symmetries
of a curvature tensor (not necessarily including the Bianchi identity)
will be called a ``vertical curvature tensor'' on the bundle.
We think of a vertical curvature tensor as prescribing the curvatures of vertical 2-planes
at zero-section.  For $p\in\Soul$ and $W,V\in E_p$ we write $k_F(W,V)=R_F(W,V,V,W)$.
As before, by parallel transporting $W$ and $V$ along geodesics from $p$ in $\Soul$, we can
think of $k_F(W,V)$ as a real valued function on a neighborhood of $p$ in $\Soul$, and
we write $\hess_{k_F(W,V)}(X)$ for the hessian of this function in the direction $X\in T_p\Soul$.

We prove that a strengthening of the necessary condition in Theorem A is
sufficient to guarantee a metric of positive curvature on the sphere bundle.

\begin{thC}
Let $\Soul$ be a compact manifold, and let $\R^k\ra E\stackrel{\pi}{\ra}\Soul$ be a vector bundle over $\Soul$.
If this bundle admits a metric $g_{\Soul}$ on $\Soul$, a Euclidean structure
$\lb\cdot,\cdot\rb$, a compatible connection $\nabla$, and a vertical curvature tensor $R_F$
such that for all $p\in\Soul$, $X,Y\in T_p\Soul$ with $X\wedge Y\neq 0$, and $W,V\in E_p$
with $W\wedge V\neq 0$:
$$
\lb (D_X\RN)(X,Y)W,V\rb^2<(|\RN(W,V)X|^2+(2/3)\hess_{k_F(W,V)}(X))\cdot k_{\Soul}(X,Y),
$$
then the unit-sphere bundle $E^1$ of $E$ admits a metric of positive curvature.
\end{thC}

Some comments about Theorem C are in order.
\begin{itemize}
\item
To prove the theorem, we construct a metric $g_E$ on $E$ for which the boundary of
a small ball about the zero-section has positive curvature.  We believe that $g_E$ can always
be constructed to be a complete metric of nonnegative curvature, but we are only able
to prove this in the special case of connection metrics, as describe in Theorem B.

\item
Guijarro and Walschap proved that if a vector bundle admits a metric of nonnegative
curvature, then so does the associated sphere bundle~\cite{GW1}.  This is because the boundary
of a small ball about the soul is convex and hence nonnegatively curved in the induced
metric.  Our theorems address the question of when this induced metric on the sphere bundle
has positive curvature.  For a metric
of nonnegative curvature on a vector bundle, the inequality in Theorem A must hold;
if in addition this inequality is strict on orthonormal vectors \{$X,Y,W,V$\},
then the induced metric on the sphere bundle must have positive curvature.

\item
The strict inequality implies
that $|\RN(W,V)X|^2\geq-(2/3)\hess_{k_F(W,V)}(X)$, with equality only when
$X=0$ or $W\wedge V=0$.  This can be thought of as a generalized fatness condition.
Because of the added generality, Derdzinski and Rigas' result does not rule out the
possibility of using Theorem C to find metrics of positive curvature on 7-dimensional
exotic spheres.

\item
In theorem A, if the metric on $M$ is such that each fiber of the projection
$\pi:M\ra\Soul$ is radially-symmetric (although not necessarily totally geodesic),
Then $(2/3)k_F(W,V)=f(p)\cdot|W\wedge V|^2$ for some function $f:\Soul\ra\R$.  In other
words, all vertical 2-planes at a fixed point $p\in\Soul$ have the same sectional curvature, so
the vertical curvature information can be described entirely by a function $f$ on $\Soul$.
In this case, the inequality of Theorem A becomes:
\begin{equation}\label{E:f}
\lb (D_X\RN)(X,Y)W,V\rb^2\leq(|\RN(W,V)X|^2+|W\wedge V|^2\cdot\hess_f(X))\cdot k_{\Soul}(X,Y).
\end{equation}
Conversely, if a vector bundle $\R^k\ra E\stackrel{\pi}{\ra}\Soul$ admits structure
$\{g_{\Soul},\lb\cdot,\cdot\rb,\nabla,f\}$ (where $\{g_{\Soul}, \lb\cdot,\cdot\rb,\nabla\}$ are
as in Theorem C, and $f:\Soul\ra\R$) such that inequality~\ref{E:f} is satisfied and is strict
for orthonormal vectors $\{X,Y,W,V\}$, then Theorem C provides a metric of positive curvature
on $E^1$ for which the fibers are round (with varying diameters).

\item
Both the inequality of Theorem A and its sharpening to a strict inequality 
for orthonormal vectors have natural interpretations.
As we will show, a mixed 2-plane $\sigma$ at the soul is a critical point of the sectional
curvature function, $\text{sec}$, on the Grassmannian $G$ of $2$-planes on $M$.
The inequality in Theorem A comes from fact that, since $M$ has nonnegative curvature,
the hessian of $\text{sec}$ at $\sigma$ must be nonnegative definite.  The strictness
of the inequality means that the
only vectors in $T_{\sigma}G$ contained in the nullspace of the hessian of $\text{sec}$
are the vectors forced to be there by Perelman's Theorem.
\end{itemize}

Our paper is organized as follows.  In sections 2 and 3 we describe the derivatives
at the soul of the $A$ and $T$ tensors of the Riemannian submersion from
an open manifold of nonnegative curvature onto its soul.  This allow us in section 4
to describe the hessian of $\text{sec}$ at a mixed 2-plane $\sigma$ at the soul. 
Theorem A is a consequence of this discussion.

In section 5 we describe how to construct a metric on a vector bundle from the data
$\{g_{\Soul},\lb\cdot,\cdot\rb,\nabla,R_F\}$ prescribed in theorem C.  Our construction
yields a ``warped connection metric'', which means a metric
formed from a connection metric by altering the fiber metrics, in our case
so that the curvatures of the fibers at the zero-section are as prescribed by $R_F$.

In section 6 we prove Theorem C by constructing a warped connection metric
on a vector bundle so that the intrinsic metric on the boundary of a small ball
about the zero-section has positive curvature.
Unfortunately, we do not know how to verify that the warped connection
metric itself has nonnegative curvature.  But in section 7 we at least show how to do this in
the case of connection metrics, thus proving Theorem B.

The author would like to thank Wolfgang Ziller and Burkhard Wilking for sharing some
insightful ideas which are incorporated into this paper.


\section{Background: The metric near the soul}
In this section, $M$ will denote an open manifold with nonnegative curvature, and
$\Soul\subset M$ will denote a soul of $M$ in the sense of~\cite{CG}.  Let $\nabla$
denote the connection in the normal bundle $\nb$ of $\Soul$ in $M$, and let $\RN$
denote it's curvature tensor.  Let $R$ denote the curvature tensor of $M$,
and denote $(E_1,E_2,E_3,E_4)=R(E_1,E_2,E_3,E_4)=\lb R(E_1,E_2)E_3,E_4\rb$, and
$k(E_1,E_2)=(E_1,E_2,E_2,E_1)$.
Define $R_{\Soul}$ and $R_F$ (respectively $k_{\Soul}$ and $k_F$)
as the restrictions of $R$ (respectively $k$) to $T\Soul$ and $\nb$.

Our proof of Theorem A relies heavily on Perelman's
resolution of the soul conjecture~\cite{P}, which states:
\begin{theorem}[Perelman]\label{T:Perelman}\hspace{.5in}
\begin{enumerate}
\item
The metric projection $\pi:M\ra\Soul$, which sends each point $p\in M$ to the point
$\pi(p)\in\Soul$ to which it is closest, is a well-defined Riemannian submersion.
\item
For any $p\in\Soul$, $X\in T_p\Soul$, and $V\in\nbp$, the surface $(s,t)\mapsto
\exp(s\cdot V(t))$ $(s,t\in\R, s\geq 0)$,
where $V(t)$ denotes the parallel transport of $V$ along the geodesic with initial
tangent vector $X$, is a flat and totally geodesic half plane (we will refer to these
surfaces as ``Perelman flats'').
\end{enumerate}
\end{theorem}

One consequence of Perelman's theorem is that mixed 2-planes at the soul are flat:

\begin{cor}\label{C:Mixed}
Let $p\in\Soul$, $X,Y\in T_p\Soul$ and $U\in\nbp$.  Then:
\begin{enumerate}
\item $R(X,U)U=R(U,X)X=0$.  In particular $k(X,U)=0$.
\item $R(X,Y)U=2R(X,U)Y$.
\end{enumerate}
\end{cor}
\begin{proof}
Part 1 was originally proved by Cheeger and Gromoll~\cite[Theorem 3.1]{CG}.  Today it is
obvious from Perelman's Theorem (at least it is obvious that $k(X,U)=0$,
but on manifold of nonnegative curvature this implies that $R(X,U)U=R(U,X)X$).
Part 2 is found in~\cite{SW}. It is a consequence of the Bianchi identity:
$$R(X,Y)U=R(X,U)Y-R(Y,U)X,$$
together with the vanishing of the mixed curvatures, which means that:
$$0=R(X+Y,U)(X+Y)=R(X,U)Y+R(Y,U)X+0+0.$$
\end{proof}

Although $\pi$ has only been proven to be $C^2$~\cite{C2},
it is clearly $C^{\infty}$ in a neighborhood of $\Soul$.
We denote by $A$ and $T$ the fundamental tensors associated to $\pi$, as defined for example
in~\cite[Chapter 9]{Besse}.  We collect in
the following lemma some facts about the $A$ and $T$ tensors at the soul.  All but part 4
of this lemma are well-known.

\begin{lem}\label{L:ATatSoul}
Let $p\in\Soul$, $X,Y\in T_p\Soul$, and $U,V,W\in\nbp$.
\begin{enumerate}
\item
$A_XY=A_XU=0.$
\item
$T_UX=T_UV=0.$
\item
$(D_VA)_XY=-\half\RN(X,Y)V$ and $(D_VA)_XU=\half\RN(V,U)X$.
\item
$(D_WT)_UV=(D_WT)_UX=0$
\end{enumerate}
\end{lem}

\begin{proof}

Part 1 is obvious.  For part 2, $T_UX=0$ because the Perelman flat through $X$ and $U$ is
totally geodesic.  Since $\lb T_UV,X\rb=-\lb T_UX,V\rb=0$, it follows that $T_UV=0$ as well.
Part 3 is an immediate consequence of~\cite[Proposition 1.7]{SW}.

For part 4, first notice that $(D_WT)_WX=0$ because the Perelman flat
through $X$ and $W$ is totally geodesic.  Next, since
$$0=\lb(D_WT)_WX,U\rb=-\lb(D_WT)_WU,X\rb=-\lb(D_WT)_UW,X\rb,$$ it follows that
$(D_WT)_UW=0$ as well.  This is a special case of part 4.  To get the general case from this
special case, we apply O'Neill's formula~\cite[9.28b]{Besse} and Corollary~\ref{C:Mixed}
as follows:  $$0=(X,W,W,U)=-\lb(D_WT)_UW,X\rb+\lb(D_UT)_WW,X\rb=\lb(D_UT)_WW,X\rb.$$
It follows from this that $(D_UT)_WW=0$.  Since $(D_UT)_{W_1}W_2$ is symmetric in
$W_1$ and $W_2$, it follows that $D_UT=0$.
\end{proof}

A different, more illuminating proof of part 4 of Lemma~\ref{L:ATatSoul} will appear
later in our proof of part 3 of Lemma~\ref{L:WarpProps}.

The following formula for the curvature of an arbitrary 2-plane at the soul
appears in~\cite[page 615]{WalSub}:
\begin{prop}[Walschap]\label{P:CurvAtSoul}
Let $p\in\Soul$, $X,Y\in T_p\Soul$, and $U,V\in\nbp$.  Then:
$$k(X+U,Y+V)  =  k_{\Soul}(X,Y) + k_F(U,V) - 3(X,Y,U,V)$$
\end{prop}

Proposition~\ref{P:CurvAtSoul} is proven by expanding the left side by linearity, and
noticing that many of the resulting terms vanish by Corollary~\ref{C:Mixed} and the fact
that the soul is totally geodesic.

\begin{cor}[Walschap]\label{C:Wal}
For all $p\in\Soul$, $X,Y\in T_p\Soul$, and $U,V\in\nbp$:
$$9(X,Y,U,V)^2\leq 4k_{\Soul}(X,Y)\cdot k_F(U,V)$$
\end{cor}
\begin{proof}
Fix $X,Y\in T_p\Soul$ and $U,V\in\nbp$.  By Proposition~\ref{P:CurvAtSoul},
$$k_{\Soul}(X,Y) + k_F(U,V) - 3(X,Y,U,V)\geq 0.$$  Of course this inequality
must remain true for any rescallings, $xX,yY,uU,vV$, of the vectors ($x,y,u,v\in\R$).
In other words, $Q(xy,uv)\geq 0$, where $Q$ is the quadratic form with matrix
$$Q=\begin{pmatrix} k_{\Soul}(X,Y)  &  -(3/2)(X,Y,U,V)\\
                    -(3/2)(X,Y,U,V) &  k_F(U,V)\end{pmatrix}$$
Hence, $Q$ is nonnegative definite and therefore has nonnegative determinate.
\end{proof}


\section{The second derivative of the $T$-tensor at the Soul}
The only new observation in the previous section was that the first derivative, $DT$,
of the $T$-tensor of $\pi$ vanishes at points of $\Soul$.
The main goal of this section is to describe the
second derivative, $D^2T$.  Since the $T$-tensor measures the failure of the
fibers to be mutually isometric, one might expect $D^2T$ at $\Soul$ to measure
the failure of the fibers to look the same at the soul; in other words, the failure
of $R_F(W_1,W_2,W_3,W_4)$ to be constant on a path in $\Soul$ along which
the sections $W_i$ of $\nb$ are parallel.  This intuition is essentially right, although
it is cleaner to describe $D^2T$ in terms of the symmetrization, $\RC_F$, of $R_F$.

We therefore begin with a discussion of symmetrization.
If $\mathbf{V}$ is a vector space with orthonormal basis $\{e_1,...,e_k\}$, and
$R$ is a curvature tensor on $\mathbf{V}$, then $\RC:S^2\mathbf{V}\ra
S^2\mathbf{V}$ commonly denotes the induces tensor on symmetric 2-forms, namely:
$$\RC(h)(U,V)=\sum_{1\leq i,j\leq k}R(e_i,U,V,e_j)\cdot h(e_i,e_j).$$
It is also useful to define
$\RC(W_1,W_2,U,V)=\RC(h)(U,V)$, where $h\in S^2\mathbf{V}$ is:
$$h(e,f)=\half(\lb e,W_1\rb\lb f,W_2\rb + \lb e,W_2\rb\lb f,W_1\rb).$$
In this way, we consider $\RC$ to be a tensor of order $4$ on $\mathbf{V}$, which has
the following simple description:
$$\RC(W_1,W_2,U,V)=\half(R(W_1,U,V,W_2)+R(W_2,U,V,W_1)).$$
The following symmetries of $\RC$ follow from the symmetries of $R$:
\begin{gather}\label{RCsymm}
\RC(W_1,W_2,U,V)=\RC(W_2,W_1,U,V)=\RC(W_1,W_2,V,U)=\RC(U,V,W_1,W_2);\\
\RC(W,W,W,U)=0.\nonumber
\end{gather}

We now give a useful description of $\RC$ when $\mathbf{V}$ is the tangent space,
$T_p M$, of a Riemannian manifold, and $R$ is the curvature tensor given by the Riemannian
metric on $M$.  First, for $W,U,V\in T_p M$, define:
\begin{equation}\label{E:Fdef}
F(W,U,V)=\lb \text{d}\exp_W U, \text{d}\exp_W V\rb.
\end{equation}
\begin{lem}\label{L:RC}
$\RC$ can be described in terms of $F$ as follows:
$$\RC(W_1,W_2,U,V)=-\left(\frac{3}{2}\right)\frac{\text{d}^2}{\text{ds}\,\text{dt}}
  \Bigr|_{s=t=0}F(tW_1+sW_2,U,V).$$
\end{lem}
\begin{proof}
Let $\Upsilon(W_1,W_2,U,V)$ denote the right side of the equation, which we wish
to prove equals $\RC(W_1,W_2,U,V)$.  Since both $\Upsilon$ and $\RC$ are symmetric
in $W_1,W_2$ and also in $U,V$, it suffices to prove that 
$\RC(W,W,U,U)=\Upsilon(W,W,U,U)$ for all $W,U\in T_p M$.
It is straightforward to see that:
$$\RC(W,W,U,U)=R(W,U,U,W)=k(W,U).$$    
It therefore remains to prove that:
$$\Upsilon(W,W,U,U)=-(2/3)k(W,U) \text{ for all vectors $W,U$}.$$
In fact, since $\Upsilon(W,W,W,U)=0$, it will suffice to verify this
when $W$ and $U$ are orthonormal, in which case $k(W,U)$ is the sectional curvature of the
2-plane $\sigma$ which they span.

Let $S$ denote the surface in $M$ obtained as the exponential image of $\sigma$.
Write the metric on $S$
in polar coordinates: $$\text{ds}^2=\text{dr}^2+f^2(r,\theta)d\theta^2,$$
where $\theta=0$ corresponds to the direction of $W$.  Let $\gamma(r)=\exp(rW)$ (in polar
coordinates, $\gamma(r)=(r,0)$). Along $\gamma$, $f$ can be expressed as:
$$f(r,0)=r\sqrt{F(rW,U,U)}.$$  The Gauss curvature
of $S$ at $(r,0)$ equals $\frac{-f_{rr}(r,0)}{f(r,0)}$, where $f_{rr}$ denotes the second
partial with respect to $r$.  The result now follows by performing the differentiation
and taking the limit as $r\ra0$.  
\end{proof}

We return to our setup where $M$ is an open manifold with nonnegative curvature,
$\Soul\subset M$ is a soul, and $R_F$ is the vertical curvature tensor.  
Let $\RC_F$ denote the fiber-wise symmetrization of $R_F$ as described above.
When the vectors $W_i\in\nbp$ are fixed, we can think of $\RC_F(W_1,W_2,W_3,W_4)$ as a real valued
function on $\Soul$ near $p$ (by parallel transporting the $W_i$ along geodesics
in $\Soul$ from $p$).  The following lemma describes $D^2T$ in terms of the gradient
of this real valued function.

\begin{lem}\label{L:DDT}
For all $p\in\Soul$ and $W_1,W_2,U,V\in\nbp$,
$$(D_{W_1}D_{W_2}T)_UV=\frac{1}{3}\text{grad}\,\RC(W_1,W_2,U,V).$$
\end{lem}

\begin{proof}
We first establish a convention for lifting vectors.  If $X\in T_p\Soul$, we denote
by $\bar{X}$ an extension of $X$ to a basic vector field on $M$.
Additionally, for $U\in\nbp$, let $\bar{U}$ be the extension $U$ to a vertical vector
field on $M$ in a neighborhood of $p$ constructed as follows.  First, extend $U$ to a vector
field along the fiber $\pi^{-1}(p)$ in a neighborhood of $p$ by defining, for each
$W\in\nbp$ with small norm, $\bar{U}|_{\exp(W)}=d\exp_W(U)$.  Then extend $U$ to a section
of $\nb$ near $p$ by parallel transporting $U$ along geodesics in $\Soul$ from $p$.
Finally, for each point $q\in\Soul$ near $p$, we extend the vector field the fiber
$\pi^{-1}(q)$ in the same way we extended it to $\pi^{-1}(p)$.

We begin by proving that for any $X\in T_p\Soul$ and $U,V,W\in\nbp$ with $|W|$ small:
\begin{equation}\label{E:TUX}
\lb T_{\bar{V}}\bar{X},\bar{U}\rb \Big|_{\exp(W)} = \half XF(W,U,V),
\end{equation}
where $F$ (which is defined in equation~\ref{E:Fdef}) is thought
of as a real-valued function on $\Soul$ near $p$ by parallel transporting $W,U,V$
along geodesics in $\Soul$ from $p$.

Notice that $\bX$ and $\bU$ commute simply because their preimages
under $d\expp$ in $T(\nb)$ commute (here $\expp:\nb\ra M$ denotes the normal exponential map).
So, letting $\bar{p}=\exp(W)$, and using the standard coordinate-free expression for the
connection, we see that:

\begin{eqnarray*}
2\lb T_{\bV}\bX,\bU\rb_{\bar{p}}
 & = & 2\lb\nabla_{\bV}\bX,\bU\rb_{\bar{p}}\\
 & = & \bX\lb\bV,\bU\rb_{\bar{p}} + \bV\lb\bU,\bX\rb_{\bar{p}} - \bU\lb\bX,\bV\rb_{\bar{p}} \\
 &   & -\lb[\bX,\bU],\bV\rb_{\bar{p}} - \lb[\bU,\bV],\bX\rb_{\bar{p}} - \lb[\bX,\bV],\bU\rb_{\bar{p}} \\
 & = & \bX\lb\bV,\bU\rb_{\bar{p}} = XF(W,U,V).
\end{eqnarray*}
This verifies equation~\ref{E:TUX}.  It follows easily that for any $W,U,V\in\nbp$ with $|W|$ small:
\begin{equation}\label{E:TUV}
(T_{\bar{V}}\bar{U})\Big|_{\exp(W)} = -\half\overline{\text{grad}F(W,U,V)}
\end{equation}
 
Finally, we use equation~\ref{E:TUV} to study $D^2T$.  To prove the lemma, it will suffice to verify
that $(D^2_WT)_UV=\frac{1}{3}\text{grad}\,\RC(W,W,U,V)$ for all $U,V,W\in\nbp$, which is done
as follows.  Let $\gamma(t)=\exp(tW)$.  Then:
\begin{eqnarray*}
\lefteqn{D_W(D_WT_UV)}\\
  & = & \Ddt\bigg|_{t=0}(D_{\bar{W}}T_{\bU}{\bV})_{\gamma(t)}-0-0\\
  & = & \Ddt\bigg|_{t=0}\left(
         \frac{\text{D}}{\text{dr}}\bigg|_{r=0}(T_{\bU}\bV)_{\gamma(t+r)}
         -(T_{(\nabla_{\bW}\bU)}\bV)_{\gamma(t)}
         -(T_{\bU}(\nabla_{\bW}\bV))_{\gamma(t)}\right)\\
  & = & \frac{\text{D}^2}{\text{dt}^2}\bigg|_{t=0}(T_{\bU}\bV)_{\gamma(t)}-0-0\\
  & = & -\half\frac{\text{D}^2}{\text{dt}^2}\bigg|_{t=0}
         (\overline{\text{grad}\,F(tW,U,V)})_{\gamma(t)}
    =   -\half\frac{\text{D}^2}{\text{dt}^2}\bigg|_{t=0}
         \text{grad}\,F(tW,U,V) \\
  & = & -\half\text{grad}\left(\frac{\text{d}^2}{\text{dt}^2}\bigg|_{t=0} F(tW,U,V)\right)
    =   \frac{1}{3}\text{grad}\,\RC(W,W,U,V).
\end{eqnarray*}
To justify the third equality above, notice that:
\begin{eqnarray*}
\Ddt\bigg|_{t=0}(T_{(\nabla_{\bW}\bU)}\bV)_{\gamma(t)}
 & = & (D_WT)_{(\nabla_W\bU)}V+T_{(\nabla_W\nabla_{\bW}\bU)}V+T_{(\nabla_W\bU)}(\nabla_W\bV)\\
 & = & 0+0+0=0
\end{eqnarray*}
\end{proof}


\section{A necessary condition for nonnegative curvature on vector bundles}
In this section we prove Theorem A by studying the derivatives of a function which
records the curvatures of a family of 2-planes.  The family begins with a mixed
2-plane at the soul and then drift so
that the base point moves away from the soul while simultaneously the 2-plane twists away
from being a mixed 2-plane.

More precisely, the set up for this section is as follows.
Let $p\in\Soul$, $X,Y\in T_p\Soul$ and
$W,V,U\in\nbp$.  Let $\gamma(t):=\exp(tW)$, and let
$X_t,Y_t,U_t,V_t$ denote the parallel transports of $X,Y,U$ and $V$ along $\gamma(t)$.
By Perelman's Theorem, $X_t$ and $Y_t$ are horizontal for all
$t\in[0,\infty)$.  In other words, parallel translation along the radial geodesic
$\gamma$ preserves the horizontal space.  It therefore must also preserve the vertical space,
so $U_t$ and $V_t$ are vertical for all $t\in[0,\infty)$.
Define:
\begin{equation}\label{E:CurvDef}
\Psi(t)=\Psi_{XYUVW}(t)=k(X_t+tU_t,tY_t+V_t)
\end{equation}
which is the unnormalized sectional curvature
of the 2-plane based at $\gamma(t)$ spanned by $X_t+tU_t$ and $tY_t+V_t$.
The special case of this construction when $U=Y=0$ was studied by Marenich in~\cite{Marenich}.
Notice that $\Psi(0)=k(X,V)=0$ by Corollary~\ref{C:Mixed}.
The goal of this section is to derive formulas for $\Psi'(0)$ and $\Psi''(0)$.
Towards this end, we write:
\begin{eqnarray}\label{tenterms}
\lefteqn{\text{\hspace{.3in}}\Psi(t) = (X_t,V_t,V_t,X_t) + t\cdot\{2(X_t,Y_t,V_t,X_t)+2(X_t,V_t,V_t,U_t)\}}\\
      & &+t^2\cdot\{(X_t,Y_t,Y_t,X_t)+(U_t,V_t,V_t,U_t)+2(X_t,Y_t,V_t,U_t)+2(X_t,V_t,Y_t,U_t)\}\nonumber\\
      & &+t^3\cdot\{2(U_t,Y_t,V_t,U_t)+2(X_t,Y_t,Y_t,U_t)\}+t^4\cdot(U_t,Y_t,Y_t,U_t)\nonumber
\end{eqnarray}

\begin{prop}\label{P:FirstDeriv}
$\Psi'(0)=0$
\end{prop}
\begin{proof}
Since $M$ has nonnegative curvature, $\Psi'(0)\geq 0$.  But if it were the case that
$\Psi'(0)>0$,
then replacing $W$ with $-W$ would yield $\Psi'(0)<0$.  Hence $\Psi'(0)=0$.

In order that our proof generalizes properly in the next section, we also
compute $\Psi'(0)$ directly.  From equation~\ref{tenterms}:
$$\Psi'(0)=2(X,Y,V,X)+2(X,V,V,U)+\dz(X_t,V_t,V_t,X_t),$$
but $(X,Y,V,X)=0$ because the soul is totally geodesic, and $(X,V,V,U)=0$ by
Corollary~\ref{C:Mixed}.  We use O'Neill's formula~\cite[Theorem 9.28c]{Besse}
to study the third term (\cite{Besse} uses a different curvature sign convention):
\begin{eqnarray*}
\lefteqn{\dz(X_t,V_t,V_t,X_t)}\\
  & = & \dz\{\lb(D_{X_t}T)_{V_t}V_t,X_t\rb
                     -\lb T_{V_t}X_t,T_{V_t}X_t\rb + \lb A_{X_t}V_t,A_{X_t}V_t\rb\}\\
  & = & \dz\lb(D_{X_t}T)_{V_t}V_t,X_t\rb - 0 + 0\\
  & = & \lb\Ddt ((D_{X_t}T)_{V_t}V_t),X\rb=\lb (D_WD_XT)_VV,X\rb=\lb (D_XD_WT)_VV,X\rb=0.
\end{eqnarray*}
\end{proof}

\begin{prop}\label{P:SecondDeriv}
\begin{eqnarray*}
\Psi''(0) & = & 2k_{\Soul}(X,Y)+2k_F(U,V)-6\lb\RN(X,Y)U,V\rb+(1/2)|\RN(W,V)X|^2\\
       &   & -2\lb(D_X\RN)(X,Y)W,V\rb +(1/3)\hess_{k_F(W,V)}(X)\\
       &   & +(4/3) D_X\RC(W,U,V,V)-(4/3)D_X\RC(W,V,U,V) 
\end{eqnarray*}
\end{prop}

\begin{proof} From equation~\ref{tenterms},
\begin{eqnarray*}
\Psi''(0) & = & 2\{(X,Y,Y,X)+(U,V,V,U)+2(X,Y,V,U)+2(X,V,Y,U)\}\\
       &   & +2\dz\{2(X_t,Y_t,V_t,X_t)+2(X_t,V_t,V_t,U_t)\}\\
       &   & +\frac{\text{d}^2}{\text{dt}^2}\bigg|_{t=0}(X_t,V_t,V_t,X_t).
\end{eqnarray*}

The top line of this expression can be simplified using part 2 of Corollary~\ref{C:Mixed}:
\begin{eqnarray*}
\lefteqn{(X,Y,Y,X)+(U,V,V,U)+2(X,Y,V,U)+2(X,V,Y,U)}\hspace{1in}\\
 & = & k_{\Soul}(X,Y)+k_F(U,V)-3(X,Y,U,V)
\end{eqnarray*}

Next, from one of O'Neill's formulas (\cite[9.28e]{Besse}):
\begin{eqnarray*}
\lefteqn{\dz(X_t,Y_t,V_t,X_t)}\\
   & = & \dz\{\lb(D_{X_t}A)_{X_t}Y_t,V_t\rb + 2\lb A_{X_t}Y_t,T_{V_t}X_t\rb\}\\
   & = & \dz\lb(D_{X_t}A)_{X_t}Y_t,V_t\rb = \lb\Ddt((D_{X_t}A)_{X_t}Y_t),V\rb\\
   & = & \lb (D_WD_XA)_XY,V\rb=\lb(D_XD_WA)_XY,V\rb\\
   & = & -\half\lb (D_X\RN)(X,Y)W,V\rb
\end{eqnarray*}

We apply another of O'Neill's formulas (\cite[9.28b]{Besse}) to simplify the next term:
\begin{eqnarray*}
\lefteqn{\dz(X_t,V_t,V_t,U_t)}\\
   & = & -\dz\{\lb(D_{V_t}T)_{U_t}V_t,X_t\rb - \lb(D_{U_t}T)_{V_t}V_t,X_t\rb\}\\
   & = & -\lb(D_WD_VT)_UV,X\rb +\lb(D_WD_UT)_VV,X\rb\\
   & = & -(1/3)\lb\grad\,\RC(W,V,U,V),X\rb + (1/3)\lb\grad\,\RC(W,U,V,V),X\rb\\
   & = & -(1/3) D_X\RC(W,V,U,V) + (1/3)D_X\RC(W,U,V,V)
\end{eqnarray*}

Finally,
\begin{eqnarray*}
\lefteqn{\frac{\text{d}^2}{\text{dt}^2}\bigg|_{t=0}(X_t,V_t,V_t,X_t)}\\
  & = & \frac{\text{d}^2}{\text{dt}^2}\bigg|_{t=0}
	 \{\lb(D_{X_t}T)_{V_t}V_t,X_t\rb - \lb T_{V_t}X_t,T_{V_t}X_t\rb + \lb
          A_{X_t}V_t,A_{X_t}V_t\rb\} \\
  & = & \lb (D^2_WD_XT)_VV,X\rb - 2\lb (D_WT)_VX,(D_WT)_VX\rb \\
  &   &  +2\lb(D_WA)_XV,(D_WA)_XV\rb\\
  & = & \lb (D_XD^2_WT)_VV,X\rb - 0 + \half|\RN(W,V)X|^2.
\end{eqnarray*}
\end{proof}

Theorem A is an immediate corollary of Proposition~\ref{P:SecondDeriv}, as we now show:

\begin{proof}[Proof of Theorem A]
Let $p\in\Soul$, $X,Y\in T_p\Soul$ and $U,V,W\in\nbp$.
Since $M$ has nonnegative curvature,
$\Psi_{XYUVW}''(t)\geq 0$.  In particular, this is true when $U=0$, which implies that the
following expression is nonnegative:
$$2k_{\Soul}(X,Y)+(1/2)|\RN(W,V)X|^2-2\lb(D_X\RN)(X,Y)W,V\rb+(1/3)\hess_{k_F(W,V)}(X).$$

Of course the same remains true for any rescallings,
$xX,yY,wW,vV$, of the vectors ($x,y,w,v\in\R)$.  In other words,
$Q(xy,xwv)\geq 0$, where $Q$ is the quadratic form with matrix:
$$
Q=\begin{pmatrix}
2k_{\Soul}(X,Y)        & \lb(D_X\RN)(X,Y)W,V\rb \\
\lb(D_X\RN)(X,Y)W,V\rb & (1/2)|\RN(W,V)X|^2+(1/3)\hess_{k_F(W,V)}(X)
\end{pmatrix}
$$
Hence, $Q$ is nonnegative definite, and it's determinate is therefore nonnegative.
This implies that:
$$\lb (D_X\RN)(X,Y)W,V\rb^2\leq(|\RN(W,V)X|^2+(2/3)\hess_{k_F(U,V)}(X))\cdot k_{\Soul}(X,Y).$$
\end{proof}

We end by mentioning a more general possible definition of $\Psi$, namely,
\begin{equation}
\Psi(t)=\Psi_{X\hat{X}YUV\hat{V}W}(t)=k(X_t+t\hat{X}_t+tU_t,tY_t+V_t+t\hat{V}_t).
\end{equation}

Although this seems more general, it is easy to see that
$\Psi_{X\hat{X}YUV\hat{V}W}'(0)=0$, and
$\Psi_{X\hat{X}YUV\hat{V}W}''(0)=\Psi_{XYUVW}''(0)$.
In other words, our derivative formulas don't notice the added generality.


\section{Warped connection metrics}
In this section we define and study a class of metrics on vector bundles called
``warped connection metrics'', which are more general than connection metrics.

Given a connection metric $g_E$ on the total space $E$ of a vector bundle,
$\R^k\ra E\stackrel{\pi}{\ra}\Soul$, we write $$TE=\Hor\oplus\V\oplus\mathbf{r}$$
for the orthogonal decomposition of the tangent bundle of $E$, where $\Hor$ is the
distribution determined by $\nabla$, $\mathbf{r}$ is the span of gradient of the distance
to the zero-section ($\mathbf{r}$ is $1$ dimensional on $E-\Soul$ and $k$ dimensional on $\Soul$),
and $\V$ describes the space of vectors tangent to the fibers of $\pi$ and
orthogonal to $\mathbf{r}$.
We make the following definition:
\begin{definition}
A warped connection metric $g_E$ on a vector bundle $\R^k\ra E\stackrel{\pi}{\ra}\Soul$
is any smooth metric obtained by starting with a connection metric and then altering the
metric arbitrarily on $\V$.
\end{definition}

For a warped connection metric, is is easy to see that $\pi$
is still a Riemannian submersion and $\exp:\nb\ra E$ is still a diffeomorphism.
Also, the zero-section, $\Soul$, is totally geodesic, and both statements of Perelman's
theorem (Theorem~\ref{T:Perelman}) are valid.
We consider the following structures on $(E,g_E)$, all defined analogous to the way they
were defined for nonnegatively curved metrics:
$g_{\Soul}$, $k_{\Soul}$, $\nabla$, $\RN$, $F$, $R_F$, $k_F$, $\RC_F$,
$A$, $T$, and $\Psi$.
For example, $F$ (which we call the warping function) is defined by the equation:
\begin{equation}\label{Fdef}
F_p(W,U,V)=\lb \text{d}\exp_WU,\text{d}\exp_WV\rb\text{ for $p\in\Soul$ and $W,U,V\in E_p$}.
\end{equation}
Notice that the vectors $\text{d}\exp_W U$ and $\text{d}\exp_W V$ are both tangent to the fibers
of $\pi$; hence, $F$ records the metrics of the fibers.  $F$ is a smooth function
from $\{(p,W,U,V)\mid p\in\Soul\text{ and }W,U,V\in\nbp\}$ to $\R$, and $F$ 
has the following properties:
\begin{enumerate}
\item $F(W,\cdot,\cdot)$ is a symmetric positive-definite bilinear form for each $W$.
\item $F(W,W,U)=F(0,W,U)=\lb W,U\rb$
\item $\frac{\text{d}}{\text{dt}}|_{t=0}F(tW,U,V)=0$
\item $\frac{\text{d}^2}{\text{dsdt}}|_{s=t=0}F(tW_1+sW_2,U,V)=
       \frac{\text{d}^2}{\text{dsdt}}|_{s=t=0} F(tU+sV,W_1,W_2)$.
\end{enumerate}

Perelman's Theorem implies that any complete metric of nonnegative curvature on a vector bundle
agrees with a warped connection metric inside of the cut-locus of the soul (this is because,
if the normal bundle, $\nb$, of the soul $\Soul$ in $M$ is endowed with its natural
connection metric, then $\exp:\nb\ra M$ preserves horizontal and vertical spaces,
and is an isometry on $\Hor$ and $\mathbf{r}$).  Guijarro proved
in~\cite{LuisTestTube} that a metric of nonnegative curvature on a vector bundle
can always be altered so that $\exp:\nb\ra M$ becomes a diffeomorphism; this altered metric is
a warped connection metric.  So, the class of warped connection metrics is general
enough to resolve Cheeger and Gromoll's question; that is, if a vector bundle admits a metric
of nonnegative curvature, then it admits a warped connection metric of nonnegative curvature.
On the other hand, the class of warped connection metrics is fairly rigid.
The next lemma says that warped connection metrics share much of the important structure of
nonnegatively curved metrics:

\begin{lem}\label{L:WarpProps}
For a warped connection metric $g_E$ on the total space $E$ of a vector bundle
$\R^k\ra E\stackrel{\pi}{\ra}\Soul$,
the following are true for all $p\in\Soul$, $X,\hat{X},Y\in T_p\Soul$ and
$W,W_1,W_2,U,V,\hat{V}\in E_p$:
\begin{enumerate}
\item
$A_p=0$ and $T_p=0$.
\item
$(D_VA)_XY=-\half\RN(X,Y)V$ and $(D_VA)_XU=\half\RN(V,U)X$.
\item
$DT_p=0$ and $(D_{W_1}D_{W_2}T)_UV=\frac{1}{3}\text{grad}\RC(W_1,W_2,U,V)$.
\item
$R(X,V)V=R(V,X)X=0$ and $R(X,Y)U=2R(X,U)Y$.
\item
$k(X+U,Y+V)=k_{\Soul}(X,Y)+k_F(U,V)-3(X,Y,U,V)$.
\item
$\Psi_{XYUVW}'(0)=\Psi_{X\hat{X}YUV\hat{V}W}'(0)=0.$
\item
$\Psi_{XYUVW}''(0)=\Psi_{X\hat{X}YUV\hat{V}W}''(0)$ is given by the equation of
Proposition~\ref{P:SecondDeriv}.
\item
The boundary of a sufficiently small ball about $\Soul$ is convex.
\end{enumerate}
\end{lem}

\begin{proof}
This lemma essentially follows from previous arguments, but one
alteration is needed.  For nonnegatively curved metrics, the fact that
$k(X,V)=0$ implies that $R(X,V)V=R(V,X)X=0$, which in turn was used to prove
that $DT_p=0$.  For general warped connection metrics,
$k(X,V)=0$, but this does not automatically imply that
$R(X,V)V=R(V,X)X=0$.  We must prove things in a different order.
First we show that $DT_p=0$.  Using O'Neill's formula, this then implies that
$R(X,V)V=R(V,X)X=0$.

To show that $DT_p=0$, let $\gamma(t)=\expp(tW)$.  Then:
\begin{eqnarray*}
D_WT_UV
  & = & \Ddt\bigg|_{t=0}(T_{\bU}{\bV})_{\gamma(t)}
    =   -\half\Ddt\bigg|_{t=0} (\overline{\text{grad}F(tW,U,V)})_{\gamma(t)}\\
  & = & -\half\Ddt\bigg|_{t=0} \text{grad}F(tW,U,V)
    =   -\half\text{grad}\left(\frac{\text{d}}{\text{dt}}\bigg|_{t=0}F(tW,U,V)\right)=0
\end{eqnarray*}

Part 8 was proven for nonnegatively curved metrics by Guijarro and Walschap
in~\cite{GW1}, and their proof remains valid for general warped connection metrics.
\end{proof}

Just as a connection metric on a vector bundle is prescribed by a Euclidean structure,
a connection, a metric on the base space, and a rotationally-symmetric metric on
$\R^k$, a warped connection metric can also be prescribed by a structures on the bundle.  
Suppose that $\R^k\ra E\stackrel{\pi}{\ra}\Soul$ is a vector bundle.  Let
$g_{\Soul}$ be a metric on $\Soul$.  Let $F$ be any smooth function from
$\{(p,W,U,V)\mid p\in\Soul\text{ and }W,U,V\in\nbp\}$ to $\R$ which has the following
two properties:
\begin{enumerate}
\item $F_p(W,\cdot,\cdot)$ is a symmetric positive-definite bilinear form for each
$p\in\Soul$ and each $W\in E_p$.
\item $F_p(W,W,U)=F_p(0,W,U)$ for each $p\in\Soul$ and each $W,U\in E_p$.
\end{enumerate}
We call $F$ a ``warping function''.  The following properties follow from the above two:
\begin{enumerate}
\item[3.] $\frac{\text{d}}{\text{dt}}|_{t=0}F(tW,U,V)=0$
\item[4.] $\frac{\text{d}^2}{\text{dsdt}}|_{s=t=0}F(tW_1+sW_2,U,V)=
       \frac{\text{d}^2}{\text{dsdt}}|_{s=t=0}F(tU+sV,W_1,W_2)$.
\end{enumerate}
To see this, notice that by property 1, $F$ induces a smooth metric on each fiber $E_p$
as follows: $\lb U,V\rb = F_p(W,U,V)$, where $U,V\in T_WE_p$, and $T_WE_p$ is identified
with $E_p$ in the obvious manner.  By property 2, the identity map from $E_p$ to $E_p$
is the exponential map with respect to this metric.  Properties 3 and 4 are now familiar facts
about metrics in polar coordinates.

$F$ determines a Euclidean structure on the bundle as follows: $\lb \cdot,\cdot\rb = F(0,\cdot,\cdot)$.
Suppose that $\nabla$ is a connection compatible with this Euclidean structure.
Then there exists a unique warped connection metric $g_E$ on $E$ for which
$\pi:(E,g_E)\ra(\Soul,g_{\Soul})$ is a Riemannian submersion with horizontal distribution
determined by $\nabla$ and with fiber metrics determined by $F$ as described above,
so that $F(W,U,V) = \lb d\exp_W U, d\exp_W V\rb$ for all $p\in\Soul$ and all
$W,U,V\in E_p$.  To construct $g_E$,
begin with the connection metric with flat fibers determined by
$g_{\Soul}$, $\lb\cdot,\cdot\rb$, and $\nabla$, and then alter the fiber metrics
according to $F$.
We call $g_E$ ``the warped connection metric on $E$ determined by the data
$\{g_{\Soul},\nabla,F\}$''.

The warping function $F$ can itself be prescribed in terms of more basic structures.
More precisely, suppose that $\lb\cdot,\cdot\rb$ is a Euclidean structure on a vector bundle,
and $R_F$ is a vertical curvature tensor on the vector bundle.  Let $\RC_F$ be the
fiber-wise symmetrization of $R_F$, as described in section 3.  Define:
$$F(W,U,V)=\lb U,V\rb - \frac{1}{3}\RC_F(W,W,U,V).$$
It is easy to verify that $F$ is smooth and has properties 1 and 2 above.
We call the resulting metric $g_E$, ``the warped connection metric on $E$ determined
by the data $\{g_{\Soul},\lb\cdot,\cdot\rb,\nabla,R_F\}$'',
even though, $g_E$ is a non-degenerate
metric only in a neighborhood of the zero-section, and not necessarily on all of $E$.  

If $R_F'$ (and $\RC_F'$) denote the vertical curvature tensor (and its symmetrization)
associated with this warped connection metric, it is clear from Lemma~\ref{L:RC}
that $\RC_F'=\RC_F$.  It follows that $R_F(W,V,V,W)=R_F'(W,V,V,W)$ for all $p\in\Soul$ and
all $W,V\in E_p$.  If $R_F$ satisfies the Bianchi identity (which we need not assume for the
previous discussion), then this implies that $R_F=R_F'$.
In other words, we succeeded in prescribing the fiber metric so that
at the zero-section its curvature tensor is $R_F$.

\section{Proof of theorem C}
\begin{proof}[Proof of Theorem C]
Suppose that $\R^k\ra E\stackrel{\pi}{\ra}\Soul$ is a
vector bundle over $\Soul$ which admits the structures
$\{g_{\Soul},\lb\cdot,\cdot\rb,\nabla, R_F\}$
so that the inequality of Theorem C is satisfied.  We wish to choose a warped connection
metric $g_E$ on $E$ such that a small sphere about $\Soul$ has positive curvature.
An obvious first try is the warped connection metric determined by the data
$\{g_{\Soul},\lb\cdot,\cdot\rb,\nabla, R_F\}$.  However, this turns out not to work.
The problem is that, since only the hessian of $R_F$ appears in the inequality of
Theorem C, this choice would provide no control over how large the vertical sectional
curvatures at the zero-section are.  For example, in the connection metric case,
the inequality is satisfied for $R_F=0$, but using flat fibers is clearly a poor choice.

We therefore modify $R_F$ to boost the sectional curvatures of vertical 2-planes
at the zero section.  When $C$ is a real number, let $R_C$ denote the vertical curvature
tensor on the vector bundle which satisfies:
$$R_C(W,U,U,W)=C\cdot|W\wedge U|^2 \text{ for all $p\in\Soul$, and $W,U\in E_p$}.$$
In other words, for each $p\in\Soul$, $(R_C)_p$ is the curvature tensor corresponding
to a point with constant sectional curvature $C$.
Let $g_E$ be the warped connection metric on $E$ which is determined by the data
$\{g_{\Soul},\lb\cdot,\cdot\rb,\nabla, R_F'\}$, where $R_F'=R_C+R_F$.
We denote by $\RC_F$  and $\RC'_F$ the symmetrizations of $R_F$ and
$R_F'$.  We denote by $k_F$ and $k_F'$ the unnormalized sectional curvatures
of $R_F$ and $R'_F$.  Notice that $\text{hess}_{k_F'(W,V)}=\text{hess}_{k_F(W,V)}$.

We will prove that for sufficiently large $C$,
the boundary of a sufficiently small ball about the zero-section of
$(E,g_E)$ has positive curvature.
\begin{claim}\label{C:1}
$C$ can be chosen sufficiently large so that the curvature of every 2-plane at every point of the zero-section $\Soul$ of $(E,g_E)$
is positive, except for the mixed 2-planes, whose curvatures are zero.
\end{claim}
\begin{proof}[Proof of Claim~\ref{C:1}]
By part 5 of Lemma~\ref{L:WarpProps} and the proof of Corollary~\ref{C:Wal},
this claim follows from the fact that $C$ can easily be chosen
so that
$$9(X,Y,U,V)^2\leq 4k_{\Soul}(X,Y)\cdot k_F'(U,V),$$
with equality only when $X\wedge Y=0$ or $U\wedge V=0$.
\end{proof}

\begin{claim}\label{C:2}
$C$ can be chosen so that if $p\in\Soul$, $X,Y\in T_p\Soul$ and $U,V,W\in E_p$ are vectors
for which $|X|=|V|=1$, $\lb W,U\rb=0$ and $\lb W,V\rb=0$, then
$$\Psi_{XYUVW}''(0)>0.$$
\end{claim}
\begin{proof}[Proof of Claim~\ref{C:2}]
By the hypothesis of Theorem C and a compactness argument, there exist $\epsilon>0$
depending only on $\{g_{\Soul},\lb\cdot,\cdot\rb,\nabla, R_F\}$ so that:
\begin{eqnarray}\label{E:epsilon}
\lefteqn{\lb (D_X\RN)(X,Y)W,V\rb^2} \\
 & & \leq(1-\epsilon)(|\RN(W,V)X|^2+(2/3)\hess_{k_F(W,V)}(X))\cdot k_{\Soul}(X,Y)\nonumber
\end{eqnarray}
for all orthonormal vectors $\{X,Y,W,V\}$ with $X,Y\in T_p\Soul$ and $W,V\in E_p$.
But since equation~\ref{E:epsilon} is invariant under
rescallings of the vectors, projection of $Y$ perpendicular to $X$, and projection of $V$
perpendicular to $W$,
this inequality is in fact valid for all (not necessarily orthonormal) vector $\{X,Y,W,V\}$.

By the argument of our proof of Theorem A, equation~\ref{E:epsilon} implies that
the following is true for all $X,Y\in T_p\Soul$ and $W,V\in E_p$:
\begin{eqnarray}\label{E:Add}
\lefteqn{$$2(1-\epsilon)k_{\Soul}(X,Y)+(1/2)|\RN(W,V)X|^2} \\
& & -2\lb(D_X\RN)(X,Y)W,V\rb+(1/3)\hess_{k_F(W,V)}(X)\geq 0.\nonumber
\end{eqnarray}
Now, part 7 of Lemma~\ref{L:WarpProps} says that $\Psi''(0)$ is given by the formula
of Proposition~\ref{P:SecondDeriv}, which we can re-write as follows:

\begin{eqnarray}\label{4lines}
\Psi''(0) & = & 2(1-\epsilon)k_{\Soul}(X,Y)+(1/2)|\RN(W,V)X|^2\\
       &   & -2\lb(D_X\RN)(X,Y)W,V\rb+(1/3)\hess_{k_F(W,V)}(X)\nonumber\\
       &   & +2\epsilon k_{\Soul}(X,Y)+2k_F(U,V)-6\lb\RN(X,Y)U,V\rb\nonumber\\
       &   & +(4/3) D_X\RC_F'(W,U,V,V)-(4/3)D_X\RC_F'(W,V,U,V)\nonumber 
\end{eqnarray}

Now let $X,Y\in T_p\Soul$ and $U,V,W\in E_p$ be vectors 
for which $|X|=|V|=1$, $\lb W,V\rb=0$ and $\lb W,U\rb=0$.  Also assume that $|W|=1$.
We must prove that $\Psi''(0)=\Psi_{XYUVW}''(0)>0$.
Equation~\ref{E:Add} says that the top two lines of equation~\ref{4lines} are nonnegative.
Further, the hypothesis of Theorem C together with a compactness argument gives that:
$$(1/2)|\RN(W,V)X|^2+(1/3)\hess_{k_F(W,V)}(X)>\delta>0,$$
for some $\delta>0$ depending only on $\{g_{\Soul},\lb\cdot,\cdot\rb,\nabla, R_F\}$.
It follows that there exist constants $\delta_1,\delta_2>0$ such that if
$|Y|\leq\delta_1$, then the sum of the terms on the first two lines of equation~\ref{4lines}
is $>\delta_2$.

Let $H$ denote the sum of the terms of the last two lines of equation~\ref{4lines}.
It will suffice to choose $C$ sufficiently large that:
\begin{enumerate}
\item if $|Y|\leq\delta_1$ then $H+\delta_2\geq 0$, and
\item if $|Y|>\delta_1$ then $H>0$.
\end{enumerate}
But notice that $D_X\RC_F'=D_X\RC_F$, and that
$$k_F(U,V)=\RC_F'(U,U,V,V)=\RC_F(U,U,V,V)+C|U\wedge V|^2.$$
From this it is straightforward to choose $C$ large enough that the above two
conditions are met.  This completes the proof of claim 2 under the added hypothesis that
$|W|=1$.  But by equation~\ref{4lines}, it is clear that
$\Psi''_{X, aY, aU, V, aW}(0)=a^2\Psi''_{X,Y,U,V,W}(0)$, which allows us to drop the 
assumption that $|W|=1$.
\end{proof}

We prove now that, if $C$ is chosen as large as required for claims 1 and 2,
then the boundary of a sufficiently small ball about $\Soul$ in $(E,g_E)$ has
positive extrinsic curvature.  By part 8 of Lemma~\ref{L:WarpProps}, it must then have
positive intrinsic curvature as well.

It is useful to consider the following manifold:
$$\Omega=\{(p,X,Y,U,V,W)\mid p\in\Soul; X,Y\in T_p\Soul; U,V,W\in E_p\}.$$
Define $f:\Omega\ra\R$ as follows:
$$f(p,X,Y,U,V,W)=k(\bar{X}+\bar{U},\bar{Y}+\bar{V}),$$
where $\{\bar{X},\bar{Y},\bar{U},\bar{V}\}$ are the lifts of $\{X,Y,U,V\}$ to $T_W E$
(via parallel translation along $W$).  Notice that $f$ has value zero on the compact submanifold
$N=\{(p,X,0,0,V,0)\in\Omega\mid |X|=|V|=1\}$.  By part 6 of Lemma~\ref{L:WarpProps},
$N$ is a critical submanifold of $f$.

Next consider the following subset of $\Omega$:
$$\Omega_{\epsilon} =
  \{(p,X,Y,U,V,W)\in\Omega\mid |W|<\epsilon, \lb W,U\rb=0\text{ and }\lb W,V\rb=0\}.$$
Notice that $N\subset\Omega_{\epsilon}$.  Even though $\Omega_{\epsilon}$ is not
a smooth submanifold of $\Omega$, we can still make the following definition.
Let $T$ denote the collection of all vectors tangent to $\Omega$ at points
of $N$ which are also ``tangent to'' $\Omega_{\epsilon}$.  In other words,
for $n\in N$ and $J\in T_n\Omega$, $J\in T$ if and only if there exists a path
$\gamma$ in $\Omega$ with $\gamma'(0)=J$ such that $\gamma([0,\delta])\subset T$ for some
$\delta>0$.

By Claim~\ref{C:2}, the Hessian of $f$ at each point of $N$ is positive-definite
in the directions of $T$.  That is, for all $J\in T$, $\text{hess}_f(J)>0$.
Since the collection of unit-directions in $T$ is compact, it follows that there
is a neighborhood of $N$ in $\Omega_{\epsilon}$ on which $f$ is strictly positive (except
on $N$ itself, where $f$ is zero).

Therefore, in the Grassmannian of 2-planes on $M$, there is a neighborhood of each
mixed 2-plane at each point of the soul in which the curvature of every 2-plane tangent to
a distance sphere about the soul is positive.
But by Claim~\ref{C:1}, every non-mixed 2-plane at the soul
has positive curvature, and therefore also has such a neighborhood.  This proves that
small distance spheres about $\Soul$ in $(E,g_E)$ have strictly positive curvature, which
completes the proof.
\end{proof}


\section{Proof of Theorem B}
In this section, we prove Theorem B.
As we mentioned in the introduction, Strake and Walschap showed by explicit computation that,
for a connection metric $g_E$ of nonnegative curvature on the total space $E$ of
a vector bundle $\R^k\ra E\stackrel{\pi}{\ra}\Soul$,
inequality~\ref{E:SWInequality} from our introduction is satisfied.
Therefore, the weaker inequality~\ref{E:ConnectInequality} is satisfied as well.
Their argument uses only the nonnegativity of 2-planes of the form $\text{span}\{X+V,Y\}$.
These 2-planes are tangent to distance spheres about $\Soul$ in $E$.
Also, their argument provides strict inequality when these 2-planes have strictly
positive curvature.  By the Gauss equation, the intrinsic curvature of such a 2-planes
(in the induced metric on the distance sphere) equals its extrinsic curvature in $(E,g_E)$.
Consequently, if some distance sphere about $\Soul$ in $E$ has strictly positive
curvature in the induced metric, then the inequality of Theorem B must be satisfied.

This proves one direction of part 1 of Theorem B.  The second direction of part 1,
as well as part 2, are proven next.  We begin with a lemma.

\begin{lem}
If a vector bundle $\R^k\ra E\stackrel{\pi}{\ra}\Soul$ admits a
warped connection metric $g_E$ which is nonnegatively curved in a neighborhood of the
zero-section $\Soul$, then it admits a complete metric $g_E'$ which has nonnegative curvature
everywhere.  If $g_E$ is a connection metric, then $g_E'$ can be chosen to be
a connection metric as well.
\end{lem}

\begin{proof}
It follows from part 8 of Lemma~\ref{L:WarpProps} that the main construction
of~\cite{LuisTestTube} can be used to modify the metric $g_E$ into a complete
metric $g_E'$ with everywhere nonnegative curvature (alternately, use the main construction
of~\cite{Kronwith}, on which~\cite{LuisTestTube} is based).
\end{proof}

\begin{proof}[Proof of Theorem B]
Let $\R^k\ra E\stackrel{\pi}{\ra}\Soul$ be a vector bundle with the structures
$\{g_{\Soul},\lb\cdot,\cdot\rb,\nabla\}$ satisfying the inequality of Theorem B.
Let $g_E$ denote a connection metric on $E$ determined by this data, whose
fiber metric is chosen so that the curvature of every vertical 2-plane at every point of
the zero-section $\Soul$ equals $C$ (for example, Walschap and Strake used the fiber metric
$dr^2+G^2(r)d\sigma^2$, where $G^2(r)=\frac{3r^2}{3+Cr^2}$).
We wish to prove that for large enough $C$, a neighborhood of $\Soul$ in $(E,g_E)$ must
have nonnegative curvature.  Together with the Lemma, this will complete the proof.

Let $p\in\Soul$, and let $W\in E_p$ be a vector with small norm.
Let $X,Y\in T_p\Soul$ and $U,V\in E_p$ with $|X|=|Y|=|U|=|V|=1$,
$\lb W,U\rb=0$ and $\lb W,V\rb = 0$.
Let $\bar{X},\bar{Y},\bar{U},\bar{V}$ denote lifts of $X,Y,U,V$ to $T_W E$
(via parallel transport along $W$).  Let $\partial_r\in T_W E$ denote the ``radial vector'',
by which we mean the unit-length vector pointing directly away from $\Soul$.
Let $\alpha,\beta,\gamma,\delta,\zeta,\eta$ be positive real numbers.  We need to insure
that $k(\alpha\bar{X}+\beta\partial_r+\gamma\bar{V},
        \delta\partial_r+\zeta\bar{Y}+\eta\bar{U})\geq 0$.
Much of the work of this calculation was done by Walschap and Strake in~\cite{Wal}.
As they showed,
\begin{eqnarray*}
\lefteqn{k(\alpha\bar{X}+\beta\partial_r+\gamma\bar{V},\delta\partial_r+\zeta\bar{Y}+\eta\bar{U})}\\
& = & k(\alpha  \bar{X}+\gamma\bar{V},\zeta\bar{Y}+\eta\bar{U})
       - \epsilon k(\alpha\bar{X},\zeta\bar{Y})
     + Q_3 (\alpha\zeta,\gamma\delta,\beta\eta),
\end{eqnarray*}
where $Q_3$ is the quadratic form with matrix:
$$
\begin{pmatrix}
\epsilon k(\bar{X},\bar{Y}) & \frac{3}{2}(\bar{X},\bar{Y},\partial_r,\bar{V}) & -\frac{3}{2}(\bar{X},\bar{Y},\partial_r,\bar{U}) \\
\frac{3}{2}(\bar{X},\bar{Y},\partial_r,\bar{V}) & k(\partial_r,\bar{V}) & (\partial_r,\bar{V},\partial_r,\bar{U}) \\
-\frac{3}{2}(X,Y,\partial_r,\bar{U}) & (\partial_r,V,\partial_r,\bar{U}) & k(\partial_r,\bar{U})
\end{pmatrix}
$$

It follows from Walschap and Strake's work that for any value of $\epsilon>0$,
the constant $C$ can be chosen so that $Q_3$ is nonnegative definite.
It therefore remains to prove that $\epsilon$ can be chosen small enough to insure that
$$k(\alpha\bar{X}+\gamma\bar{V},\zeta\bar{Y}+\eta\bar{U})
 -\epsilon k(\alpha\bar{X},\zeta\bar{Y})\geq 0.$$
We will modify the argument in Theorem C by which we proved that
small distance spheres are positively curved to get this little bit more which is required.

We can use the value of $\epsilon$ from equation~\ref{E:Add}, which simplifies in the 
connection metric case to:
$$2(1-\epsilon)k_{\Soul}(X,Y)+(1/2)|\RN(W,V)X|^2-2\lb(D_X\RN)(X,Y)W,V\rb\geq 0.$$
Modify the proof of Theorem C by defining $f:\Omega\ra\R$ as follows:
$$f(p,X,Y,U,V,W)=k(\bar{X}+\bar{U},\bar{Y}+\bar{V})-\epsilon k(\bar{X},\bar{Y}).$$
It is still easy to see that $N$ is a critical submanifold of $f$,
and that the hessian of $f$ is positive definite in directions of $T$.
The proves that, in the Grassmannian of 2-planes on $M$, there is a neighborhood
of each mixed 2-plane at each point of $\Soul$ on which the sectional curvature
function is nonnegative.  Since non-mixed 2-planes at points of $\Soul$
also have such neighborhoods,
this proves that a neighborhood of $\Soul$ in $(E,g_E)$ has nonnegative curvature, which
completes the proof.
\end{proof}
\bibliographystyle{amsplain}

\end{document}